\documentclass[a4paper,10pt]{article}
\usepackage{authblk}

\usepackage[ngerman,english]{babel}
\usepackage[latin1]{inputenc}
\usepackage{csquotes}
\usepackage[normalem]{ulem}
\usepackage{amsfonts,amsmath,amsthm}
\usepackage{empheq}
\usepackage{cancel}
\usepackage[titletoc,title]{appendix}

\usepackage{cases}
\usepackage{mathabx}
\usepackage{xfrac}
\usepackage{fancyhdr}
\usepackage{color}
\usepackage[colorlinks]{hyperref}
\definecolor{blue75}{rgb}{0,0,.75}
\definecolor{green75}{rgb}{0,.75,0}
\hypersetup{colorlinks=true, urlcolor=blue75,linkcolor=blue75,citecolor=green75,pdfstartview=FitB,bookmarksopen=true,bookmarksopenlevel=1}
\usepackage[a4paper, left=1.8cm, right=1.8cm, top=1.8cm,bottom=1.8cm]{geometry}
\usepackage{constants}
\newcommand{\parenthezises}[1]{\arabic{#1}}
\newconstantfamily{C}{
symbol=C,
format=\parenthezises,
reset={section}
}
\newconstantfamily{M}{
symbol=M,
format=\parenthezises,
reset={section}
}
\newconstantfamily{a}{
symbol=a,
format=\parenthezises,
reset={section}
}
\usepackage{enumerate}

\usepackage{graphicx}
\graphicspath{ {images/} }
\usepackage{wrapfig}
\usepackage{figbib}
\allowdisplaybreaks
\usepackage[capitalise]{cleveref}

\crefdefaultlabelformat{{\it #2#1#3}}

\crefname{equation}{}{}

\crefname{enumi}{}{}
\creflabelformat{enumi}{{(#2#1#3)}}

\crefname{section}{{\it Section}}{{\it Sections}}
\crefname{subsection}{{\it Subsection}}{{\it Subsections}}
\crefname{subsubsection}{{\it Paragraph}}{{\it Paragraphs}}
\crefname{table}{{\it Table}}{{\it Tables}}
\crefname{figure}{{\it Figure}}{{\it Figures}}

\crefname{Theorem}{{\it Theorem}}{{\it Theorems}}

\crefname{Definition}{{\it Definition}}{{\it Definitions}}

\crefname{Lemma}{{\it Lemma}}{{\it Lemmas}}

\crefname{Assumptions}{{\it Assumptions}}{{\it Assumptions}}

\theoremstyle{definition}

\crefname{Remark}{{\it Remark}}{{\it Remarks}}

\crefname{Notation}{{\it Notation}}{{\it Notations}}

\crefname{Example}{{\it Example}}{{\it Examples}}

\usepackage{titling}
\begin{document}

\newcommand{\FOO}{n(n+1)}
\definecolor{darkgreen}{rgb}{0.1,0.7,0.2}

\definecolor{orange}{rgb}{1,0.5,0}

\newcommand{\cb}{\color{blue}}
\newcommand{\cred}{\color{red}}
\newcommand{\cmg}{\color{magenta}}
\newcommand{\cg}{\color{darkgreen}}
\newcommand{\cy}{\color{cyan}}
\newcommand{\co}{\color{orange}}
\newcommand{\D}{\mathbb{D}}
\newcommand{\R}{\mathbb{R}}
\newcommand{\N}{\mathbb{N}}
\newcommand{\F}{\mathbb{F}}
\newcommand{\Om}{\Omega }
\newcommand{\eps}{\epsilon }
\newcommand{\J}{\mathbb{J}}

\def\diam{\operatorname{diam}}
\def\dist{\operatorname{dist}}
\def\ess{\operatorname{ess}}
\def\inner{\operatorname{int}}
\def\osc{\operatorname{osc}}
\def\sign{\operatorname{sign}}
\def\supp{\operatorname{supp}}
\newcommand{\BMO}{BMO(\Omega)}
\newcommand{\LOne}{L^{1}(\Omega)}
\newcommand{\LOnen}{(L^{1}(\Omega))^n}
\newcommand{\LTwo}{L^{2}(\Omega)}
\newcommand{\Lq}{L^{q}(\Omega)}
\newcommand{\Lp}{L^{p}(\Omega)}
\newcommand{\Lpn}{(L^{p}(\Omega))^n}
\newcommand{\LInf}{L^{\infty}(\Omega)}
\newcommand{\HOneO}{H^{1,0}(\Omega)}
\newcommand{\HTwoO}{H^{2,0}(\Omega)}
\newcommand{\HOne}{H^{1}(\Omega)}
\newcommand{\HTwo}{H^{2}(\Omega)}
\newcommand{\HmOne}{H^{-1}(\Omega)}
\newcommand{\HmTwo}{H^{-2}(\Omega)}

\newcommand{\LlogL}{L\log L(\Omega)}

\def\avint{\mathop{\,\rlap{-}\!\!\int}\nolimits}

\numberwithin{equation}{section}

\title{Mathematical models for cell migration: a nonlocal perspective}

\makeatletter
\let\@fnsymbol\@alph
\makeatother

\author{
   Li Chen\thanks{Mathematisches Institut, Universität Mannheim, A5 6, 68131 Mannheim, Germany, \href{mailto:chen@math.uni-mannheim.de}{chen@math.uni-mannheim.de}}, 
   \quad Kevin Painter\thanks{Department of Mathematics \& Maxwell Institute, Heriot-Watt University,
   Edinburgh, Scotland, EH14 4AS, \href{mailto:K.Painter@hw.ac.uk}{K.Painter@hw.ac.uk}},\quad Christina Surulescu\thanks{Felix-Klein-Zentrum für Mathematik, Technische Universität Kaiserslautern, Paul-Ehrlich-Str. 31, 67663 Kaiserslautern, Germany, \href{mailto:surulescu@mathematik.uni-kl.de}{surulescu@mathematik.uni-kl.de}
   },\quad and Anna Zhigun\thanks{School of Mathematics and Physics, Queen's University Belfast, University Road, Belfast BT7 1NN, Northern Ireland, UK, \href{mailto:A.Zhigun@qub.ac.uk}{A.Zhigun@qub.ac.uk}
   }
}
\date{\vspace{-5ex}}
\setlength{\droptitle}{-4em} 
\maketitle

\begin{abstract}
We provide a review of recent advancements in nonlocal continuous models for migration, mainly from the perspective of its involvement in embryonal development and cancer invasion. Particular emphasis is placed on spatial nonlocality occurring in advection terms, used to characterise a cell's motility bias according to its interactions with other cellular and acellular components in its vicinity (e.g., cell-cell and cell-tissue adhesions, nonlocal chemotaxis), but we also shortly address spatially nonlocal source terms. Following a brief introduction and motivation, we give a systematic classification of available PDE models with respect to the type of featured   nonlocalities and review some of the mathematical challenges arising from such models, with a focus on analytical aspects.	
\end{abstract}

\section{Introduction}\label{intro}
Collective movement arises when individuals correlate their motion 
with others, generating migration at a population level. Paradigms include 
flocks, swarms and crowds \cite{sumpter2010}, but it also occurs for 
bacteria \cite{ariel2013}, embryonic populations,
immune and invading cancer cells \cite{friedl2012,mayor2016}.  
Scales span enormous ranges, from a few cells clustered over a few microns to 
millions or billions of organisms distributed over kilometres, e.g. large-scale fish 
schools \cite{makris2009} and locust swarms \cite{rainey1967}.
An adoption of theoretical approaches has helped understand these phenomena. 
{\em Agent-based modelling} (ABMs) is a popular approach, with its
individual-level representation facilitating data fitting. For cell populations, 
ABM approaches range greatly in sophistication, including single or multi-site 
cellular automata \cite{graner1992,deutsch2005}, descriptions of cells as 
overlapping spheres \cite{drasdo2005}, deformable ellipsoids \cite{palsson2000} or 
dynamic boundaries \cite{schaller2005}. For organisms, collective
movement models are often founded on point-based individuals moving with velocities 
determined by their interactions with neighbours (see the review in \cite{berdahl2018}).

\noindent
Despite their many advantages, problems persist with ABMs that motivate 
complementary approaches. First, a lack of standard analytical methods leads to
heavy reliance on computation which, inevitably, becomes
burdensome as population size increases. Second, how should one compare the results 
emerging following different approaches applied to the same problem, 
e.g. between a lattice- and off-lattice model used to describe cell sorting behaviour? 
Precise quantitative matching is clearly unrealistic, so when can one state that two methods 
generate equivalent behaviour? Third, different implementations of the same method 
can also generate quantitatively distinct results when applied to the same problem 
\cite{rejniak2007b}.  This typically escalates with the sophistication/detail of the 
ABM, with variations arising from, say, ambiguously stated assumptions 
or distinctions in the numerical implementation. Overall, these issues highlight 
the  general challenge of appropriately ``benchmarking'' ABMs,  and we refer to \cite{osborne2017} for a more detailed consideration. 

\noindent
 While it would be disingenuous to state that {\em continuous models} are 
free from such issues, in principal their solutions are reproducible: 
well-posed problems generate unique solutions for a given set of initial 
conditions. Further, with their roots in classical theory, well-developed 
analytical methods exist that provide generic insights: the analysis 
necessary to demonstrate the self-organising capacity of Turing's 
counterintuitive reaction and diffusion theory 
of morphogenesis \cite{turing1952} is not restricted to precise reactions, 
parameters, etc. Phenomenological derivations start with a mass conservation 
equation, where movement is modelled via stipulating an appropriate {\em flux}. 
Coupled to reaction/birth/death processes, governing equations are stated for the 
key variables (cells, organisms, chemicals,  etc.), each represented by 
continuous density distributions. Models derived in this way typically fall into the
class of  reaction-diffusion-advection (RDA) equations,
\begin{equation}\label{adr}
\dfrac{\partial u}{\partial t} = \overbracket{\nabla \cdot \left( D(\cdot) \nabla u \right)}^{\mbox{Diffusion}} - \overbracket{\nabla \cdot \left({\bf {a}} (\cdot) u\right)}^{\mbox{Advection}} + \overbracket{f(u,\cdot)}^{\mbox{Reaction}} \,,
\end{equation}
where $u(x,t)$ denotes the population density at position $x$ at time $t$. Diffusion 
describes a non-oriented dispersal process, for example due to simple 
random meandering by individuals and is characterised by diffusion coefficient $D(\cdot)$.
Advection could be passive (e.g. environmental flow) or 
due to active navigation by individuals, and is described by an advective velocity 
${\bf {a}} (\cdot)$. Reaction describes the population birth/death, etc. A vast 
number of models fall in the above class, including numerous 
landmark works: textbooks such as \cite{okubo2001,murray2} offer scope of this framework.
Models of type (\ref{adr}) can also be derived as the  
continuous limiting equation of a biased random walk description for 
biological particle movement (e.g. see \cite{codling2008}).

\noindent
Models in the RDA framework typically have a {\em local} nature, i.e. terms that depend 
pointwise. For example, in the well-known Keller-Segel model \cite{keller1970} for chemotaxis, 
the  advection term describes population drift along a local chemoattractant gradient: 
specifically, ${\bf {a}} = \chi \nabla v$, for some chemoattractant $v$ and function 
$\chi$. Effectively, cells (or animals) are assumed to detect and migrate in the direction of 
a local 
gradient. This is often logical, viewed at a macroscopic level: cells such as leukocytes
orient according to the concentration difference of attractant across its body axis, but 
at the scale of a tissue this can be regarded as a pointwise calculation.

\noindent
Local assumptions may not, however, always hold or be convenient. Population densities may be 
high: classical diffusive fluxes (e.g. Fickian) assume diluteness, and at high densities 
the impact of long range effects may be important \cite{murray2}. Moreover, many particles 
sense the environment over extended regions: filopodia/cytonemes permit cells to 
detect signals multiple cell diameters away \cite{kornberg2014}; sensory organs grant 
organisms with highly nonlocal perception fields (e.g.  \cite{payne1971,lima1996,fagan2017,johnston2019}). Approximating information 
originating over large regions to, say, a local gradient, could clearly be overly-reductive. 
Dispersal distances may also be nonlocal, for example seeds can be transported significant 
distances from source while various studies have implicated ``L\'evy-type'' behaviour in 
migration paths, where short range movements are interspersed with occasional long 
transits (e.g. \cite{harris2012}). Local formulations can also create analytical problems, 
exemplified in the ``blow-up'' phenomena in certain formulations of chemotaxis models 
(e.g. see \cite{BBTW}). Here, the coupling between a population's pointwise production 
of its own attractant and movement up the local gradient leads to runaway aggregation 
and singularity formulation. Such phenomena are powerful indicators of the 
inherent self-organisation, yet formation of infinite cell densities is, ultimately, 
unrealistic. 

\noindent
These considerations and others have led to a range of spatially-nonlocal RDA 
models, and their modelling  and mathematical properties have attracted 
significant interest. This brief survey focuses on some aspects of modelling 
through such a framework. Nonlocality is, of course, a broad concept and can be 
included in various ways, for example into any or all of the diffusion, advection or 
reaction terms. We primarily focus on the use of {\em nonlocal advection} models
that feature {\em spatial integral operators} inside advection terms. These
have typically been developed to replace the gradient-type terms often used to 
describe taxis-type movement and, in particular, have come into vogue as a 
method of modelling collective movement processes in cells and organisms.

\section{Applications in development and cancer}\label{SecDev}

Nonlocal advection models have received considerable attention for their capacity to include
cell-cell (and cell-matrix) adhesion into models for tissue dynamics. Adhesion occurs 
when juxtaposing membranes link certain transmembrane adhesion 
proteins, fastening cells together and forming clusters \cite{alberts2015}. Moreover, 
cell-cell adhesion confers self-organisation, with famous studies revealing how 
mixed cell types can self-rearrange into distinct configurations, 
implying a capacity to ``recognise'' others of same type \cite{townes1955}. The Differential 
Adhesion Hypothesis (DAH) of Steinberg \cite{steinberg1963} suggested that distinct adhesion 
can provide this ``tissue-affinity'', with the ratio of self- to cross-adhesion 
strengths determining the configuration; various experiments corroborate this theory
(e.g. \cite{steinberg2007}). 

\noindent
Models of adhesion should ideally exhibit clustering/sorting, and many ABMs indeed 
reproduce these phenomena (e.g. see \cite{osborne2017}). The discrete cell 
representation is optimal: adhesion easily enters as an attracting force 
over a range of cell-cell separations, coalescing cells until their
compression generates a counteracting repulsion. Incorporating adhesion into continuous 
models, however, can prove challenging. Attempts starting from an 
initial discrete random walk process have certainly generated continuous 
models, yet these can be ill-posed (backward diffusion) or seemingly 
incapable of displaying more complicated behaviour such as sorting
(e.g. \cite{anguige2008,johnston2012,johnston2013}). 

\noindent
Phenomenological approaches founded on nonlocal concepts appear to be more 
successful. Such models capture cell-neighbour 
interactions through the proposed movement of cells according to the 
density of others in their vicinity. An early model of this type was proposed 
in \cite{sekimura1999}, although subsequent analysis focused on a localised 
form derived under expansion. The nonlocal model for adhesion proposed in 
\cite{APS06} was explored regarding its ability to recapitulate the sorting 
behaviour predicted by the DAH and its relative success has led to various extensions: 
\cite{gerisch2010} performed a more comprehensive analysis;  
\cite{15MuTo,carrillo2019} replaced the overly-reductive linear diffusion 
terms with nonlinear forms, generating the sharp cell boundaries 
often observed experimentally; \cite{painter2015} extended to more general 
cell-cell contact phenomena, for example allowing repulsive 
interactions as found in Eph-Ephrin interactions \cite{taylor2017}; the model 
of \cite{ko2019} has extended to allow dynamic adhesion regulation.

\noindent
Typical applications lie in morphogenesis and cancer. The former has 
witnessed nonlocal advection models used to describe somitogenesis \cite{armstrong2009}, 
mesenchymal condensation in early limb development \cite{glimm2014,bhat2019}, neuronal 
positioning in early brain development \cite{matsunaga2017,trush2019} and 
zebrafish gastrulation \cite{ko2019}. Notably, many of these studies integrate 
modelling with experimental data. The formulation of nonlocal advection models for 
cancer invasion has addressed the question of how cell-cell and cell-matrix adhesion interact 
with other mechanisms to facilitate cancer invasion, 
e.g. \cite{gerisch2008,kim2009,PAS10,andasari2011,DTGC14,ESuSt2016,bitsouni2017}.
As one example, the study of \cite{DTGC14} recapitulates various observed 
tumour infiltrative patterns, as well as the characteristic 
morphologies of ductal carcinomas and fibroadenomas. Other cellular applications of 
nonlocal advection models include the interactions between liver hepatocyte 
and stellate cells for {\em in vitro} culture systems \cite{green2010}. Nonlocal models of 
cell migration and spread including adhesion, have also been extended to 
account for further structure, such as cellular age and the level of bound receptors, see
\cite{DTGC14,DTGC17,dyson-webb,ESuSt2016}. Including variables characterising 
subcellular dynamics opens the way for multiscality. 

\noindent
Nonlocal advection models have also been applied extensively to 
problems of animal movement, particularly animal swarming/flocking behaviour. The 
pioneering model of \cite{mogilner1999} featured a nonlocal advection based on a 
convolution, modelling the attracting and repelling interactions between neighbouring 
swarm members. This model has sparked various extensions and significant analysis, for 
example see \cite{Lee_etal01,topaz2006,bernoff2013,fetecau2011,fetecau2013,evers2017}. 
In the context of swarming, hyperbolic approaches have been developed in 
which nonlocal interactions are included in the turning behaviour of swarm members, 
allowing extensions to orientation alignment (see the review in \cite{eftimie2012}). 
Nonlocal advection models have also been used to incorporate perceptual range into the model \cite{fagan2017,johnston2019}, i.e. animal movement according to information drawn 
from potentially large regions of their environment.

\section{Classes of nonlocal models for cell migration}\label{SecClasses}
We can extend (\ref{adr}) to a general RDA equation of the form \eqref{EqLoc}, describing the evolution of a subpopulation density $u_i$ as a part of an ensemble  $u=(u_1,\dots,u_n)$ of $n\in\N$ components representing cell densities, densities of a surrounding fibrous environment (e.g., natural or artificial tissue), concentrations of nutrients and chemical signals, etc.:
\begin{align}
 \partial_t u_i=\nabla\cdot\left(a_{i0}(u)\nabla u_i\right)-\nabla\cdot\left(\sum_{j=1}^{m-1} a_{ij}(u)\nabla b_{ij}(u)\right)+a_{im}(u).\label{EqLoc}
\end{align}
Here $\nabla=\nabla_x$ is the spatial gradient, $m\in\N$, and the coefficients have the following meaning: $a_{i0}(u)$ is the diffusion coefficient (normally nonnegative), $a_{ij}(u)$ and $b_{ij}(u)$ for $j\in\{1,\dots,m-1\}$ describe  tactic sensitivities and signal functions, respectively, and, finally, $a_{im}(u)$ is the reaction-interaction term. As previously remarked, nonlocality can be introduced in multiple ways into such PDEs. Often, it takes the form of an integral operator w.r.t. time $t$  and/or position $x$ in a spatial set $O\subset\R^d$, $d\in\N$, but other  independent variables (e.g., orientation/speed or age/phenotype/individual state, etc.) can also be involved. A typical spatial nonlocal operator can be described as follows:
\begin{align}
 {\cal I}v(x):=\int_{O}J(x,y)v(y)\,dy,\nonumber
\end{align}
where $J$ is some kernel defined in $O\times O$. If, for instance, $O=\R^d$ and $J=J(x-y)$, then the so-called convolution notation is used:   
\begin{align}
 {\cal I}v=J\star v.\nonumber
\end{align}
It can be seen e.g., as the combined ability (over the whole spatial region $O$) of some extracellular trait (mediated by a density distribution function $J$) and some quantity $v$ (density/volume fraction/etc.), to determine the cell density at a specific location $x$.

\noindent
Nonlocalities of orders zero, one, or two can be distinguished according to whether a coefficient function, a first-, or a second-order differential operator is replaced by a nonlocal one. For example, a zero-order nonlocality is present if an $a_{ij}$ is made dependent upon ${\cal I} u$. Moreover,  
nonlocality can be introduced into the reaction, taxis, or diffusion terms, leading to another 
possible  classification. In the subsequent text we address these and other possibilities in 
more detail.

\subsection{Spatial nonlocality in  advection terms}\label{subsec:nonloc-drift}
There are (at least) four ways to include a nonlocality  into the  advective flux, see \cref{Tab1}.
\begin{table}[h]
\begin{tabular}{|r|c|c|}
\hline
 Integral operator & Examples& {References}\\\hline
 { is placed before $\nabla$}&{ $(J\star v_1)\nabla v_2$}&  \cite{Lee_etal01}\\\hline
is placed inside $\nabla$&
 {$\nabla(J\star v)$}& {\cite{MCOe05}}\\\hline 
 replaces  $\nabla$ & ${\cal A}_{r}v(x)=\frac{1}{r} \avint_{B_r}v(x+\xi)\frac{\xi}{|\xi|}F_r(|\xi|)\,d\xi$&adhesion velocity {\cite{APS06,AAE05}}\\
 &$\mathring{\nabla}_{r}v(x)=\frac{n}{r} \avint_{S_r}v(x+\xi)\xi\,d_{{\cal S}_r}$&nonlocal {chemo}taxis \cite{HilPSchm,othmer-hillen2}\\\hline
 is {applied to} $\nabla$&{${\cal T}_r \nabla v(x)=\frac{1}{r}\int_0^1\avint_{B_{r}} (\nabla v(x+s\xi )\cdot \xi )\frac{\xi }{|\xi |}F_r(|\xi|)\,d\xi ds$}&\cite{KSZ}\\
 &{${\cal S}_r \nabla v(x)=\frac{n}{r}\int_0^1\avint_{S_{r}} (\nabla v(x+s\xi )\cdot \xi )\xi \,d_{{{\cal S}_r}} ds$}&\cite{KSZ}\\\hline
\end{tabular}         
\caption{Nonlocal modifications of the gradient operator  applied to a function $v$ (or $v=(v_1,v_2)$).}\label{Tab1}
\end{table} 
Hereafter $B_r$ and $S_r$ denote the open  $d$-dimensional ball and the $(d-1)$-dimensional sphere, respectively, which are centred at the origin and have radius $r$, termed the sensing radius.  The operator $\avint$ denotes the usual averaging over the set upon which the integration takes place. For the precise mathematical formulations consult the references in Table \ref{Tab1}.
Constructions in lines 1 and 2 in the Table can be viewed as zero order nonlocalities.  The former describes, e.g., the situation of long-range interactions of individuals having density $v_1$ with their environment containing a signal of concentration $v_2$ (think of cells extending protrusions towards sites with higher concentrations of some chemoattractant, i.e. directing themselves towards the gradient of such concentrations). If the chemical signal itself is assumed to move much faster than the cells -which is often the case-,  then $v_2$ can actually be expressed as a function of $v_1$, possibly in a nonlocal way, too, thus leading to a flux of the form $(J_1\star v_1)\nabla (J_2\star  v_1)$, as e.g., in \cite{Lee_etal01}. This corresponds to direct, long-range intraspecific interactions. Line 2 in Table \ref{Tab1} refers e.g., to the case of individuals (cells, ants,...) moving in a collective way, thereby perceiving and correspondingly adapting to regions with large crowd density.\footnote{Thereby $J\star v$ can be seen to represent the average density felt by the individual.} Concerning the remaining lines of Table \ref{Tab1},  an operator ${\cal M}\in \{{\cal A}_{r},\mathring{\nabla}_{r},{\cal T}_r \nabla,{\cal S}_r \nabla\}$ can be used to include a nonlocality of  first order. 
A basic model example of the latter case is given by a system of two equations
\begin{subequations}\label{nltaxis}
	\begin{align}
	\partial_t u_1=&\nabla\cdot\left(a_{10}(u)\nabla u_1-{ a_{11}(u)} {\cal M}(b_{11}(u))\right)+a_{12}(u),\label{nlu1}\\
	\partial_t u_2=&a_{20}\Delta u_2+a_{21}(u),\label{nlu2}
	\end{align}
\end{subequations}
{equipped with suitable initial and boundary conditions.} 
It can describe growth and motility of a single cell population {of density} $u_1$ biased by {intra-} and interspecies interaction{s}  and/or a signal {concentration} $u_2$. The latter is either a diffusing chemoattractant/-repellent if $a_{20}>0$, or, if $a_{20}=0$, an insoluble cue - usually {a non-diffusing} polymeric matrix such as tissue fibers.  Further components can be {included into the system, e.g., other cell populations and other soluble/insoluble signals. 

\noindent
A nonlocal chemotaxis model was introduced in \cite{othmer-hillen2} and further studied in \cite{Butten,HilPSchm,HP09,loy}. { Such settings} can be derived from position- or velocity jump processes under adequate assumptions, e.g., constant {$r$} for shrinking spatial mesh size or nonlocal sensing introducing a bias of higher order w.r.t. $r$. {This leads to the operator $\mathring{\nabla}_{r}$  in the  advection  term.} {Cell-cell and/or cell-tissue interactions are usually characterised by a so-called adhesion operator ${\cal A}_{r}$ involving a suitable function $F_r$. The latter represents the distance-dependent magnitude of the interaction force.}
We refer to \cite{APS06,Butten,loy} and references therein for formal deductions of {such models. Other versions} characterising the nonlocal space-time dynamics of one or several interacting species (cell populations, soluble and insoluble signals) have also been addressed  \cite{APS06,gerisch2008,PAS10,SGAP09}. 

\noindent
Very recently, a model class was introduced \cite{KSZ}, which uses ${\cal T}_r \nabla$ {(resp. ${\cal S}_r \nabla$)} rather than ${\cal A}_{r}$ {(resp. $\mathring{\nabla}_{r}$)}.  There, it was  pointed out that on the one hand
\begin{align}
 &{\cal A}_ru={\cal T}_r (\nabla u),\qquad \mathring{\nabla}_{r} u={\cal S}_r (\nabla u)\qquad\text{in }\Omega_r{:=\{x\in\Omega\ :\  \dist(x,\partial\Omega)>r\}},\nonumber
\end{align}
whereas, on the  other hand, e.g. for $\Omega=(-1,1)$ and $u\equiv 1$ in $\Omega$
\begin{align}
 {\cal T}_{r}(u')\equiv0 \equiv u',\qquad \int_{-1}^1|{\cal A}_ru|\,dx=1\qquad\text{for }r\in(0,1).\label{ExAT}
\end{align}
 In \cite{KSZ}  $\Om_r$ was termed  \textit{domain of restricted sensing} since there cells a priori cannot directly  perceive signals from outside the domain of interest $\Omega$. For $r\to 0$  it tends to cover the whole of $\Omega$. In contrast, a cell inside the  $r$-thick boundary layer $\Om \setminus \Om _r$ can potentially reach beyond $\partial \Om$. Of course, if $r$ is larger than the diameter of $\Omega$, then each cell  can do that. 
However, if the population is kept in a Petri dish or it is  confined within comparatively hard barriers, e.g. bone material, then  the cell flux through the boundary $\partial\Omega$  vanishes. This leads to cell densities such as $u$ from above. As \cref{ExAT} shows, in such cases the outputs under operators ${\cal A}_r$ and ${\cal T}_r \nabla$ are equal in $\Om _r$, but may disagree substantially inside $\Om \setminus \Om _r$ even for very small $r$.  In case of impenetrable boundaries and $r$ close to zero, the study in \cite{KSZ} supports the idea that cells 
actively adjust their movement after suitably sampling signal gradients rather than densities. We refer to that reference for a detailed discussion.

\noindent
Other continuous models have been obtained by starting from a particle description, e.g. accounting for long range attraction and short range repulsion between individuals in a population alongside Brownian dispersal. In the limit of sufficiently large populations these lead to nonlinear PDEs for one-component models \cite{MCOe05,Burger-capasso07} containing, for instance, a degenerate diffusion $a_{10}(u)=u$ as well as an operator ${\cal J}$ in the  advection.
Further models in this category have been proposed 
in \cite{BDF08,burger-haskovec-wolfram}. Models accounting for cell interactions with attraction or repulsion have also been studied in \cite{carrillo2019,15MuTo}. A related approach \cite{MFG14} 
employs an off-lattice ABM and derives a continuum approximation able to account for correlations between moving cells. A mean-field approximation of the evolution equations obtained for one- and two-cell density functions  {starting from} Langevin equations for cell movement leads to a PDE akin to the more common adhesion models  {from} above. Models with similar mathematical structure are also used to describe crowd dynamics, flocking, or  swarming, often referred to as self-organisation models, see \cite{bellomo-dogbe,eftimie,muntean-toschi} and references therein.

\subsection{Further types of  spatial or other nonlocality}\label{subsec:other-nonloc} 


Replacing the usual Laplace operator in a diffusion term with a fractional Laplacian (see e.g. \cite{friedman2012}) is another way of including spatial nonlocality within motility terms and exemplifies a second order nonlocality. Such models account for dispersal of individuals performing L\'evy flights rather than Brownian motion,  cf. 
\cref{intro}. Systems describing competition between locally and nonlocally dispersing populations were developed and studied in \cite{Valdinoci2017, kao2010}.


\noindent
Nonlocalities introduced into reaction-interaction terms can still affect cell motion, albeit indirectly. Indeed, cell proliferation and decay (alongside intra- and interspecific interactions) lead to local changes in densities, which flows into the density-dependent coefficients. From a modelling perspective, this accounts for population pressure,  competition for resources, cooperation in signal transmission, differentiation, and/or tissue degradation, etc. But even when motion coefficients do not depend on the population density, local vs. nonlocal source terms may lead to different overall evolution, see the discussion below. In a broader framework, classical reaction terms in population dynamics have been introduced  in \cite{fisher1937,kpp} and they are local. For the emergence and evolution of a single biological species the typical  choice is 
 \begin{align}
a(u)= \mu u^\alpha(1-u)-\gamma u.        \label{localR}                                                                                                                                                                                                         \end{align}
For $\alpha = 1$, 
growth is proportional to the population density and limited by competition for available resources.  The case $\alpha >1$ accounts for advantages of clustering together or organizing in groups. This applies to cells \cite{wang-collective}, but also to  sexual reproduction (case $\alpha =2$) or swarming of animals. 

\noindent
Individuals, of course, typically perceive information related to occupancy, biochemical cues, etc. within a neighborhood centred on their current position. Thus, local terms like \cref{localR} have been recently replaced by nonlocal ones. The best known example of the resulting equation is 
\begin{equation}\label{NonlocalKPP}
\partial_t u=\Delta u+\mu u^\alpha\left(1-J\star u^\beta\right)-\gamma u,
\end{equation}
where $J$ is a kernel {as} 
above, and $\alpha, \beta, \mu, \gamma $ are constants. 
Here the nonlocality is of  order zero.
Similar reaction terms have been used, e.g., to describe natural selection of tumour cells leading to the emergence of therapy-resistant clones  \cite{Lorz:2011hl,Lorz:2013hq}. 
Further examples of nonlocal source terms, not necessarily connected to biological applications, are of the form ${a}(u)= f(u) + I(g(u))$, where $I(\zeta ):=\int _\mathcal O\zeta (y,t)dy$, see e.g. 
\cite{Pao92,dlx03,bds93}.  In a biological context such terms can account for both local and nonlocal interactions between cells and their surroundings. We refer to \cite{kav-suzuki,vol2} for a rich variety of nonlocal reaction models in engineering and biology.

\noindent
Several model classes have also been developed featuring integral terms that  describe 
nonlocality with respect to one or several other variables, including age, phenotype, 
internal cell state, velocity, etc. They include the large class of structured population 
models \cite{magal-ruan}, as well as kinetic transport equations (KTE) (and in particular the 
so-called kinetic theory of active particles framework, see \cite{BBGO17,BB19} and references therein). Under appropriate conditions, models with spatial nonlocality can be (formally) derived from 
KTEs, see  e.g., \cite{LCS-nonlocal,othmer-hillen2}. \\  

\section{Local vs. nonlocal models: Mathematical aspects}\label{sec:challenge}

This section briefly discusses relevant qualitative results. We focus on analysis pertaining to just two model classes:  equations featuring nonlocal reaction and local diffusion, e.g. \cref{NonlocalKPP}, and  settings that involve a first order nonlocality to model a process such as adhesion or nonlocal chemotaxis  (cf. \cref{subsec:nonloc-drift}). Our motivation for this focus is as follows. On the one 
hand, the most straightforward way of accounting for nonlocality is via a zeroth order in the 
source terms.\footnote{In this review we use the standard designations 'reaction' and 'source' for all terms containing zero order derivatives. 
} Thus, understanding and overcoming challenges met when analysing such equations is essential for developing a general mathematical theory applicable to nonlocal problems.  
Consequently, the basic representative of the class, equation \cref{NonlocalKPP}, has 
received significant attention by analysts. While models involving first-order nonlocalities 
have received considerably less study, they are particularly relevant for applications,
particularly in the context of collective motion phenomena (cf. \cref{SecDev}).

\subsection{Analysis of models with spatial nonlocalities in reaction terms}\label{subsec:nonloc-react}
The analysis of reaction-diffusion equations featuring nonlocalities in source terms is highly challenging, in large part down to classical techniques that rely on comparison principles being no longer valid. A general theory seems presently out of reach, since the analysis heavily
depends on the exact form of involved nonlocality, where key features of the corresponding 
settings are revealed, for example see \cite{Pao92,dlx03,bds93}. If one includes a parameter where, 
as it is formally sent to zero, the nonlocal equation becomes local, then one can expect that 
results for the local equation can be suitably generalised to the nonlocal setting. 
As for the corresponding local case, studies of general nonlocal models such as \eqref{NonlocalKPP} include results on global well-posedness, blow-up, and stationary solutions. Specific solutions, such as stationary, radially symmetric, travelling wave solutions, or monotone wave fronts have also received attention due to their relevance in applications. 

\noindent
To exemplify, consider the relatively well understood nonlocal Fisher-KPP equation \cref{NonlocalKPP} for the case $\gamma=0$. For $J\equiv1$, which corresponds to the situation of blind competition, and
with general $\alpha,\beta\geq1$, a global bounded solution has been shown to exist both for 
bounded and unbounded domains \cite{BCL,BC}. When the kernel $J$ is replaced by the Dirac delta function, \cref{NonlocalKPP} reduces to a classical, local reaction-diffusion equation. There, results on 
global well-posedness, asymptotic stability of nontrivial stationary solutions, as well as other 
solution behaviours such as hair trigger effect\footnote{meaning that  an initially very small cell density can evolve in the long time into a cell mass completely filling the space, i.e. at maximum density}, extinction and quenching, have been intensively investigated, see for example \cite{fisher1937,ref3,Lou}. If instead $J>0$ in a ball of positive radius, then the nonlocality can have a profound impact. For instance, the constant solution $u\equiv 1$ can lose the stability of the corresponding local case with a periodic-in-space stationary solution bifurcating from it \cite{ref13,ref19,ref21}. This phenomenon has been observed in the study of travelling wave solutions, and numerically tested for the time dependent version in \cite{LCS-nonlocal}. On the other hand, if $J$ has an everywhere-positive Fourier transform or if 
it approximates the Dirac delta function, then there are travelling waves connecting $u=0$ and $u=1$ for $\alpha=1$ (see \cite{BNPR}), and \cite{LCS-nonlocal} shows that for $1\leq \alpha<1+\frac{2\beta}{N}$ the hair trigger effect appears, while for large $\mu$ values $u=1$ can indeed become unstable and Turing patterns  occur \cite{NPT}. Similar results have been obtained for the bistable case \cite{LC}. As observed in \cite{LCS-nonlocal}, the concrete solution behaviour, in particular w.r.t. pattern formation, depends strongly on the shape of the interaction kernel.
Even for \cref{NonlocalKPP} the integral kernel must be fixed to study in detail long-time behaviour. 
For systems of PDEs with nonlocalities in the reaction terms the situation is even more complicated and, to our knowledge, there has been no breakthrough in the study of behaviour in this context.


\subsection{Analysis of models with spatial nonlocalities in advection terms}
\noindent { The} rigorous analysis of local RDA systems has enjoyed great popularity over recent decades. The Keller-Segel systems are among the best studied \cite{BBTW,Horstmann,ZhigunKS}, a model class corresponding to  ${\mathcal M}=\nabla$ in \cref{nltaxis}.  
In contrast, only a few studies consider problems including one of the four nonlocal operators introduced in \cref{Tab1} that lead to first order nonlocalities. 
At a general level, combining local diffusion with nonlocal advection appears to preclude the 
existence of an energy functional satisfying a precise dissipation identity, as known for 
various formulations of local Keller-Segel model and providing a key for their {analysis}. 
Owing to this drawback, only settings where the nonlocal advection is effectively dominated by 
diffusion have been {investigated} so far. This is generally the case when { the} operators ${\cal A}_{r}$ or $\mathring{\nabla}_{r}$ are involved, since they replace a differential operator by an integral one, leading to an increase (rather than a decrease) in regularity. In the absence of other effects this allows well-posedness to be established. Moreover, it turns out that the uniform boundedness of 
solutions can be { guaranteed} under quite general assumptions, including even cases where the corresponding local system exhibits finite time blow-up. Even situations in which $a_{10}-a_{11}\partial_{u_1}b_{11}$ is  { somewhere negative} can be { covered.} In the corresponding 
local setting this implies negative self-diffusion and, generally, non-existence of solutions. 
A detailed analysis of a {nonlocal} chemotaxis system was carried out in  \cite{HilPSchm}. Several studies, in particular \cite{dyson10,dyson13,HillButten2019,SGAP09,ESuSt2016,AAE05}, address 
equations or systems {featuring} the adhesion operator ${\cal A}_r$ or its extension to a possibly unbounded sensing region \cite{WHP16}. Some works exploit specific solutions which are particularly  relevant for applications, including steady states and their stability, existence of travelling waves, etc., see  \cite{Xiang} and \cite{HilPSchm,dyson10,dyson13,Butten,OuZhang} for models with $\mathring{\nabla}_{r}$ and ${\cal A}_r$, respectively. 

\noindent
Overall, operators $\mathring{\nabla}_{r}$ and ${\cal A}_r$ form a powerful alternative to the local gradient, particularly as they allow modelling a broader range of aggregative mechanisms without fear
of potential blow-up. Moreover, as formal Taylor expansions performed in \cite{Hillen2007} and \cite{gerisch2008} respectively indicate, $\mathring{\nabla}_{r}$ and ${\cal A}_r$ approach the
local gradient $\nabla u$ for some fixed smooth  $u$ and vanishing $r$. 
In \cite{HilPSchm} 
the question was therefore raised concerning convergence of solutions to a family of nonlocal  
chemotaxis systems as $r\rightarrow0$. This corresponds to the sensing region of a cell almost shrinking to its respective position, i.e. the sensing is  effectively local.
However, as the example from \cref{subsec:nonloc-react} indicates, blow-up may appear in the gradient limit on the boundary of the spatial domain. 
Using ${\cal T}_r \nabla$ or  ${\cal S}_r \nabla$ instead excludes this undesired effect. 
These operators are, however, computed based on the gradient and they are closer to it both 
quantitatively and qualitatively. Consequently, the domination of diffusion over advection 
demands much stronger conditions on coefficients $a_{ij}$ and $b_{ij}$. Suitable conditions 
have been found and existence and rigorous convergence (of a subsequence) of solutions 
proved in  \cite{KSZ}.

\noindent
{The issue of connecting spatially local and nonlocal models acting on the same (macroscopic) scale has also been addressed, e.g., in \cite{Lee_etal01} (upon performing an adequate scaling) and, as mentioned above, in \cite{Hillen2007,gerisch2008} upon Taylor approximations (for small $r$) of functions inside the nonlocal operators. Those deductions are, however, formal, whereas \cite{KSZ} provides a rigorous approach. }





\section{Outlook}
Several challenges arise in connection with nonlocal models, some of which we already mentioned. Here we focus on models for cell migration, but most mathematical issues also apply to systems of this type characterising other real-world phenomena.

\noindent
From the \textit{modelling} viewpoint, the settings can be extended to account for various aspects of cell migration and growth. For instance, tumor heterogeneity can be w.r.t. cell phenotypes, motility, treatment response,  etc.; each of these are  influenced by the composition of the tumor microenviroment, which in turn is dynamically modified by the cells, according to their population behavior. This results  in  ODE-PDE systems with intricate couplings and nonlinearities, even if only spatial nonlocality is considered. Including several populations of cells structured by further variables, as addressed at the end of \cref{subsec:other-nonloc}, leads to multiscale descriptions, involving hyperbolic and/or parabolic PDEs with various nonlocalities. The latter can also occur in a pure macroscopic framework with only spatial nonlocality. When the cell densities evolve in a bounded domain one has to provide adequate boundary conditions. Depending on the complexity of the system accounting for interactions of cells between themselves and their surroundings, deriving them together with the population level dynamics is often nontrivial and calls for a careful modelling starting on lower scales and performing appropriate upscalings. Connections between local and nonlocal settings retain their relevance also in this context.  From the \textit{numerical} viewpoint, nonlocal models present significant challenges: integrating across a nonlocal region carries a substantial extra burden over classical local RDA models, compounded as one moves into higher (e.g. three) dimensions. Numerically efficient techniques can be developed (e.g. see \cite{gerisch2008,gerisch2010}), although they typically rely on, e.g. convenient boundary conditions or static sensing regions. Continued development of efficient methodologies is therefore a must for further, more intricate applications.

\noindent
From an \textit{analytical} viewpoint, it is desirable to support initially formal  deductions by performing a rigorous limit procedure wherever it is possible. Notwithstanding, qualitative properties, such as the well-posedness, the long-time behaviour including the possibility of a blow-up, the limit behaviour w.r.t. to some vanishing parameter, etc.,  need to be addressed for the  resulting models. Overall, these key aspects have remained open for many  cell migration models, and that includes even local, single-scale ones.  Introducing a nonlocality  into a well-understood local model can lead to additional  challenges since it breaks the original structure, see  the discussions in  \cref{sec:challenge}.
%

\phantomsection

\begin{thebibliography}{100}

\bibitem{AAE05}
M.~Adioui, O.~Arino, and N.~El~Saadi.
\newblock A nonlocal model of phytoplankton aggregation.
\newblock {\em Nonlinear Anal. Real World Appl.}, 6(4):593--607, 2005.

\bibitem{alberts2015}
B.~Alberts, A.~Johnson, J.~Lewis, D.~Morgan, M.~Raff, K.~Roberts, and
  P.~Walter.
\newblock {\em Molecular Biology of the Cell}.
\newblock Gardland Science, 6th edition, 2015.

\bibitem{andasari2011}
V.~Andasari, A.~Gerisch, G.~Lolas, A.~P. South, and M.~A. Chaplain.
\newblock Mathematical modeling of cancer cell invasion of tissue: biological
  insight from mathematical analysis and computational simulation.
\newblock {\em J. Math. Biol.}, 63(1):141--171, 2011.

\bibitem{anguige2008}
K.~Anguige and C.~Schmeiser.
\newblock {{A} one-dimensional model of cell diffusion and aggregation,
  incorporating volume filling and cell-to-cell adhesion}.
\newblock {\em J. Math. Biol.}, 2008.

\bibitem{ariel2013}
G.~Ariel, A.~Shklarsh, O.~Kalisman, C.~Ingham, and E.~Ben-Jacob.
\newblock From organized internal traffic to collective navigation of bacterial
  swarms.
\newblock {\em New J. Phys.}, 15(12):125019, 2013.

\bibitem{APS06}
N.~J. Armstrong, K.~J. Painter, and J.~A. Sherratt.
\newblock A continuum approach to modelling cell-cell adhesion.
\newblock {\em J. Theoret. Biol.}, 243(1):98--113, 2006.

\bibitem{armstrong2009}
N.~J. Armstrong, K.~J. Painter, and J.~A. Sherratt.
\newblock Adding adhesion to a chemical signaling model for somite formation.
\newblock {\em Bull. Math. Biol.}, 71(1):1--24, 2009.

\bibitem{ref3}
D.~Aronson and H.~Weinberger.
\newblock Multidimensional nonlinear diffusion arising in population genetics.
\newblock {\em Adv. Math.}, 30(1):33 -- 76, 1978.

\bibitem{BBGO17}
N.~Bellomo, A.~Bellouquid, L.~Gibelli, and N.~Outada.
\newblock {\em A quest towards a mathematical theory of living systems}.
\newblock Modeling and Simulation in Science, Engineering and Technology.
  Birkh\"{a}user/Springer, Cham, 2017.

\bibitem{BBTW}
N.~Bellomo, A.~Bellouquid, Y.~Tao, and M.~Winkler.
\newblock Toward a mathematical theory of {K}eller-{S}egel models of pattern
  formation in biological tissues.
\newblock {\em Math. Models Methods Appl. Sci.}, 25(9):1663--1763, 2015.

\bibitem{BB19}
N.~Bellomo and F.~Brezzi.
\newblock Towards a multiscale vision of active particles.
\newblock {\em Math. Models Methods Appl. Sci.}, 29(4):581--588, 2019.

\bibitem{bellomo-dogbe}
N.~Bellomo and C.~Dogb\'e.
\newblock On the modeling of traffic and crowds: A survey of models,
  speculations, and perspectives.
\newblock {\em SIAM Rev.}, 53(3):409--463, 2011.

\bibitem{berdahl2018}
A.~M. Berdahl, A.~B. Kao, A.~Flack, P.~A.~H. Westley, E.~A. Codling, I.~D.
  Couzin, A.~I. Dell, and D.~Biro.
\newblock Collective animal navigation and migratory culture: from theoretical
  models to empirical evidence.
\newblock {\em Philos. Trans. Roy. Soc. London Ser. B}, 373(1746):20170009,
  2018.

\bibitem{BNPR}
H.~Berestycki, G.~Nadin, B.~Perthame, and L.~Ryzhik.
\newblock The non-local {F}isher-{KPP} equation: travelling waves and steady
  states.
\newblock {\em Nonlinearity}, 22(12):2813--2844, 2009.

\bibitem{bernoff2013}
A.~J. Bernoff and C.~M. Topaz.
\newblock Nonlocal aggregation models: A primer of swarm equilibria.
\newblock {\em SIAM Rev.}, 55(4):709--747, 2013.

\bibitem{bhat2019}
R.~Bhat, T.~Glimm, M.~Linde-Medina, C.~Cui, and S.~A. Newman.
\newblock Synchronization of {H}es1 oscillations coordinates and refines
  condensation formation and patterning of the avian limb skeleton.
\newblock {\em Mech. Devel.}, 156:41--54, 2019.

\bibitem{BC}
S.~Bian and L.~Chen.
\newblock A nonlocal reaction diffusion equation and its relation with {F}ujita
  exponent.
\newblock {\em J. Math. Anal. Appl.}, 444(2):1479--1489, 2016.

\bibitem{BCL}
S.~Bian, L.~Chen, and E.~A. Latos.
\newblock Global existence and asymptotic behavior of solutions to a nonlocal
  {F}isher-{KPP} type problem.
\newblock {\em Nonlinear Anal.}, 149:165--176, 2017.

\bibitem{bitsouni2017}
V.~Bitsouni, M.~A.~J. Chaplain, and R.~Eftimie.
\newblock Mathematical modelling of cancer invasion: the multiple roles of
  tgf-$\beta$ pathway on tumour proliferation and cell adhesion.
\newblock {\em Math. Models Methods Appl. Sci.}, 27(10):1929--1962, 2017.

\bibitem{ref13}
N.~F. Britton.
\newblock Spatial structures and periodic travelling waves in an
  integro-differential reaction-diffusion population model.
\newblock {\em SIAM J. Appl. Math.}, 50(6):1663--1688, 1990.

\bibitem{bds93}
C.~Budd, B.~Dold, and A.~Stuart.
\newblock Blowup in a partial differential equation with conserved first
  integral.
\newblock {\em SIAM J. Appl. Math.}, 53(3):718--742, 1993.

\bibitem{Burger-capasso07}
M.~Burger, V.~Capasso, and D.~Morale.
\newblock On an aggregation model with long and short range interactions.
\newblock {\em Nonlinear Anal. Real World Appl.}, 8(3):939 -- 958, 2007.

\bibitem{BDF08}
M.~Burger and M.~{Di Francesco}.
\newblock Large time behavior of nonlocal aggregation models with nonlinear
  diffusion.
\newblock {\em Netw. Heterog. Media}, 3(4):749--785, 12 2008.

\bibitem{burger-haskovec-wolfram}
M.~Burger, J.~Ha\v{s}kovec, and M.-T. Wolfram.
\newblock Individual based and mean-field modeling of direct aggregation.
\newblock {\em Physica D}, 260:145 -- 158, 2013.
\newblock Emergent Behaviour in Multi-particle Systems with Non-local
  Interactions.

\bibitem{HillButten2019}
A.~Buttensch\"{o}n and T.~Hillen.
\newblock Nonlocal adhesion models for microorganisms on bounded domains,
  preprint.

\bibitem{Butten}
A.~Buttensch\"{o}n, T.~Hillen, A.~Gerisch, and K.~J. Painter.
\newblock A {\it space-jump} derivation for non-local models of cell-cell
  adhesion and non-local chemotaxis.
\newblock {\em J. Math. Biol.}, 76(1-2):429--456, 2018.

\bibitem{carrillo2019}
J.~A. Carrillo, H.~Murakawa, M.~Sato, H.~Togashi, and O.~Trush.
\newblock A population dynamics model of cell-cell adhesion incorporating
  population pressure and density saturation.
\newblock {\em J. Theor. Biol.}, 474:14--24, 2019.

\bibitem{codling2008}
E.~A. Codling, M.~J. Plank, and S.~Benhamou.
\newblock Random walk models in biology.
\newblock {\em J. Roy. Soc. Interface}, 5:813--834, 2008.

\bibitem{dlx03}
W.~Deng, Y.~Li, and C.~Xie.
\newblock Semilinear reaction-diffusion systems with nonlocal sources.
\newblock {\em Math. Comput. Modelling}, 37(9-10):937--943, 2003.

\bibitem{deutsch2005}
A.~Deutsch and S.~Dormann.
\newblock {\em Cellular automaton modeling of biological pattern formation}.
\newblock Springer, 2005.

\bibitem{DTGC14}
P.~Domschke, D.~Trucu, A.~Gerisch, and M.~A.~J. Chaplain.
\newblock Mathematical modelling of cancer invasion: implications of cell
  adhesion variability for tumour infiltrative growth patterns.
\newblock {\em J. Theoret. Biol.}, 361:41--60, 2014.

\bibitem{DTGC17}
P.~Domschke, D.~Trucu, A.~Gerisch, and M.~A.~J. Chaplain.
\newblock Structured models of cell migration incorporating molecular binding
  processes.
\newblock {\em J. Math. Biol.}, 75(6-7):1517--1561, 2017.

\bibitem{drasdo2005}
D.~Drasdo and S.~H{\"o}hme.
\newblock A single-cell-based model of tumor growth in vitro: monolayers and
  spheroids.
\newblock {\em Phys. Biol.}, 2(3):133, 2005.

\bibitem{dyson13}
J.~Dyson, S.~Gourley, and G.~Webb.
\newblock A non-local evolution equation model of cell-cell adhesion in higher
  dimensional space.
\newblock {\em J. Biol. Dyn.}, 7:68--87, 2013.

\bibitem{dyson10}
J.~Dyson, S.~A. Gourley, R.~Villella-Bressan, and G.~F. Webb.
\newblock Existence and asymptotic properties of solutions of a nonlocal
  evolution equation modeling cell-cell adhesion.
\newblock {\em SIAM J. Math. Anal.}, 42(4):1784--1804, 2010.

\bibitem{dyson-webb}
J.~Dyson and G.~F. Webb.
\newblock A cell population model structured by cell age incorporating
  cell-cell adhesion.
\newblock In {\em Mathematical oncology 2013}, Model. Simul. Sci. Eng.
  Technol., pages 109--149. Birkh\"{a}user/Springer, New York, 2014.

\bibitem{eftimie2012}
R.~Eftimie.
\newblock Hyperbolic and kinetic models for self-organized biological
  aggregations and movement: a brief review.
\newblock {\em J. Math. Biol.}, 65(1):35--75, 2012.

\bibitem{eftimie}
R.~Eftimie.
\newblock {\em Hyperbolic and Kinetic Models for Self-organised Biological
  Aggregations. A Modelling and Pattern Formation Approach}.
\newblock Springer, Cham, 2018.

\bibitem{ESuSt2016}
C.~Engwer, C.~Stinner, and C.~Surulescu.
\newblock On a structured multiscale model for acid-mediated tumor invasion:
  the effects of adhesion and proliferation.
\newblock {\em Math. Models Methods Appl. Sci.}, 27(7):1355--1390, 2017.

\bibitem{evers2017}
J.~H. Evers, R.~C. Fetecau, and T.~Kolokolnikov.
\newblock Equilibria for an aggregation model with two species.
\newblock {\em SIAM J. Appl. Dyn. Syst.}, 16(4):2287--2338, 2017.

\bibitem{fagan2017}
W.~F. Fagan, E.~Gurarie, S.~Bewick, A.~Howard, R.~S. Cantrell, and C.~Cosner.
\newblock Perceptual ranges, information gathering, and foraging success in
  dynamic landscapes.
\newblock {\em Amer. Nat.}, 189:474--489, 2017.

\bibitem{fetecau2013}
R.~C. Fetecau and Y.~Huang.
\newblock Equilibria of biological aggregations with nonlocal
  repulsive--attractive interactions.
\newblock {\em Physica D}, 260:49--64, 2013.

\bibitem{fetecau2011}
R.~C. Fetecau, Y.~Huang, and T.~Kolokolnikov.
\newblock Swarm dynamics and equilibria for a nonlocal aggregation model.
\newblock {\em Nonlinearity}, 24(10):2681, 2011.

\bibitem{fisher1937}
R.~A. Fisher.
\newblock The wave of advance of advantageous genes.
\newblock {\em Ann. Eugen.}, 7(4):355--369, 1937.

\bibitem{friedl2012}
P.~Friedl, J.~Locker, E.~Sahai, and J.~E. Segall.
\newblock Classifying collective cancer cell invasion.
\newblock {\em Nat. Cell Biol.}, 14(8):777, 2012.

\bibitem{friedman2012}
A.~Friedman.
\newblock P{DE} problems arising in mathematical biology.
\newblock {\em Netw. Heterog. Media}, 7(4):691--703, 2012.

\bibitem{ref19}
S.~Genieys, V.~Volpert, and P.~Auger.
\newblock Pattern and waves for a model in population dynamics with nonlocal
  consumption of resources.
\newblock {\em Math. Model. Nat. Phenom.}, 1(1):65--82, 2006.

\bibitem{gerisch2008}
A.~Gerisch and M.~A.~J. Chaplain.
\newblock Mathematical modelling of cancer cell invasion of tissue: local and
  non-local models and the effect of adhesion.
\newblock {\em J. Theoret. Biol.}, 250(4):684--704, 2008.

\bibitem{gerisch2010}
A.~Gerisch and K.~J. Painter.
\newblock Mathematical modelling of cell adhesion and its applications to
  developmental biology and cancer invasion.
\newblock {\em Cell mechanics: from single scale-based models to multiscale
  modeling}, 2:319--350, 2010.

\bibitem{glimm2014}
T.~Glimm, R.~Bhat, and S.~A. Newman.
\newblock Modeling the morphodynamic galectin patterning network of the
  developing avian limb skeleton.
\newblock {\em J. Theor. Biol.}, 346:86--108, 2014.

\bibitem{ref21}
S.~A. Gourley.
\newblock Travelling front solutions of a nonlocal {F}isher equation.
\newblock {\em J. Math. Biol.}, 41(3):272--284, 2000.

\bibitem{graner1992}
F.~Graner and J.~A. Glazier.
\newblock Simulation of biological cell sorting using a two-dimensional
  extended potts model.
\newblock {\em Phys. Rev. Lett.}, 69(13):2013, 1992.

\bibitem{green2010}
J.~E.~F. Green, S.~L. Waters, J.~P. Whiteley, L.~Edelstein-Keshet, K.~M.
  Shakesheff, and H.~M. Byrne.
\newblock Non-local models for the formation of hepatocyte--stellate cell
  aggregates.
\newblock {\em J. Theor. Biol.}, 267(1):106--120, 2010.

\bibitem{harris2012}
T.~H. Harris, E.~J. Banigan, D.~A. Christian, C.~Konradt, E.~D.~T. Wojno,
  K.~Norose, E.~H. Wilson, B.~John, W.~Weninger, A.~D. Luster, et~al.
\newblock Generalized l{\'e}vy walks and the role of chemokines in migration of
  effector cd8+ t cells.
\newblock {\em Nature}, 486(7404):545, 2012.

\bibitem{Hillen2007}
T.~Hillen.
\newblock A classification of spikes and plateaus.
\newblock {\em SIAM Rev.}, 49(1):35--51, 2007.

\bibitem{HilPSchm}
T.~Hillen, K.~Painter, and C.~Schmeiser.
\newblock Global existence for chemotaxis with finite sampling radius.
\newblock {\em Discrete Contin. Dyn. Syst. Ser. B}, 7(1):125--144.

\bibitem{WHP16}
T.~Hillen, K.~Painter, and M.~Winkler.
\newblock Global solvability and explicit bounds for non-local adhesion models.
\newblock {\em Eur. J. Appl. Math.}, 29(04):1--40, 2018.

\bibitem{HP09}
T.~Hillen and K.~J. Painter.
\newblock A user's guide to {PDE} models for chemotaxis.
\newblock {\em J. Math. Biol.}, 58(1-2):183--217, 2009.

\bibitem{Horstmann}
D.~Horstmann.
\newblock From 1970 until present: the {K}eller-{S}egel model in chemotaxis and
  its consequences. {I}.
\newblock {\em Jahresber. Deutsch. Math.-Verein.}, 105(3):103--165, 2003.

\bibitem{johnston2019}
S.~T. Johnston and K.~J. Painter.
\newblock The impact of short-and long-range perception on population
  movements.
\newblock {\em J. Theor. Biol.}, 460:227--242, 2019.

\bibitem{johnston2012}
S.~T. Johnston, M.~J. Simpson, and R.~E. Baker.
\newblock Mean-field descriptions of collective migration with strong adhesion.
\newblock {\em Physical Review E}, 85(5):051922, 2012.

\bibitem{johnston2013}
S.~T. Johnston, M.~J. Simpson, and M.~J. Plank.
\newblock Lattice-free descriptions of collective motion with crowding and
  adhesion.
\newblock {\em Physical Review E}, 88(6):062720, 2013.

\bibitem{kao2010}
C.-Y. Kao, Y.~Lou, and W.~Shen.
\newblock Random dispersal vs. non-local dispersal.
\newblock {\em Discrete Contin. Dyn. Syst.}, 26(2):551--596, 2010.

\bibitem{kav-suzuki}
N.~I. Kavallaris and T.~Suzuki.
\newblock {\em Non-local partial differential equations for engineering and
  biology}, volume~31 of {\em Mathematics for Industry (Tokyo)}.
\newblock Springer, Cham, 2018.
\newblock Mathematical modeling and analysis.

\bibitem{keller1970}
E.~Keller and L.~Segel.
\newblock Initiation of slime mold aggregation viewed as an instability.
\newblock {\em J.\ Theor.\ Biol.}, 26:399--415, 1970.

\bibitem{kim2009}
Y.~Kim, S.~Lawler, M.~O. Nowicki, E.~A. Chiocca, and A.~Friedman.
\newblock A mathematical model for pattern formation of glioma cells outside
  the tumor spheroid core.
\newblock {\em J. Theor. Biol.}, 260(3):359--371, 2009.

\bibitem{ko2019}
J.~M. Ko and D.~Lobo.
\newblock Continuous dynamic modeling of regulated cell adhesion.
\newblock {\em bioRxiv}, 2019.

\bibitem{kpp}
A.~Kolmogorov, I.~Petrovsky, and N.~Piskunov.
\newblock Investigation of the equation of diffusion combined with increasing
  of the substance and its application to a biology problem.
\newblock {\em Bull. Moscow State Univ. Ser A: Math and Mech}, 1:1--25, 1937.

\bibitem{kornberg2014}
T.~B. Kornberg and S.~Roy.
\newblock Cytonemes as specialized signaling filopodia.
\newblock {\em Development}, 141:729--736, 2014.

\bibitem{KSZ}
M.~Krasnianski, C.~Surulescu, and A.~Zhigun.
\newblock Nonlocal and local models for taxis in cell migration: a rigorous
  limit procedure, 2019.

\bibitem{Lee_etal01}
C.~Lee, M.~Hoopes, J.~Diehl, W.~Gilliland, G.~Huxel, E.~Leaver, K.~McCann,
  J.~Umbanhovar, and A.~Mogilner.
\newblock Non-local concepts and models in biology.
\newblock {\em J. Theor. Biol.}, 210(2):201 -- 219, 2001.

\bibitem{LCS-nonlocal}
J.~Li, L.~Chen, and C.~Surulescu.
\newblock Global existence, asymptotic behavior, and pattern formation driven
  by the parametrization of a nonlocal fisher-kpp problem.

\bibitem{LC}
J.~Li, E.~Latos, and L.~Chen.
\newblock Wavefronts for a nonlinear nonlocal bistable reaction-diffusion
  equation in population dynamics.
\newblock {\em J. Differential Equations}, 263(10):6427--6455, 2017.

\bibitem{lima1996}
S.~L. Lima and P.~A. Zollner.
\newblock Towards a behavioral ecology of ecological landscapes.
\newblock {\em Trends Ecol. \& Evol.}, 11:131--135, 1996.

\bibitem{Lorz:2013hq}
A.~Lorz, T.~Lorenzi, M.~E. Hochberg, J.~Clairambault, and B.~Perthame.
\newblock Populational adaptive evolution, chemotherapeutic resistance and
  multiple anti-cancer therapies.
\newblock {\em ESAIM Math. Model. Numer. Anal.}, 47(2):377--399, 2013.

\bibitem{Lorz:2011hl}
A.~Lorz, S.~Mirrahimi, and B.~Perthame.
\newblock Dirac mass dynamics in multidimensional nonlocal parabolic equations.
\newblock {\em Comm. Partial Differential Equations}, 36(6):1071--1098, 2011.

\bibitem{Lou}
Y.~Lou, T.~Nagylaki, and W.-M. Ni.
\newblock An introduction to migration-selection {PDE} models.
\newblock {\em Discrete Contin. Dyn. Syst. A}, 33(10):4349--4373, 2013.

\bibitem{loy}
N.~Loy and L.~Preziosi.
\newblock Kinetic models with non-local sensing determining cell polarization
  and speed according to independent cues.
\newblock {\em Journal of Mathematical Biology}, 2019.

\bibitem{magal-ruan}
P.~Magal and S.~Ruan.
\newblock {\em Structured Population Models in Biology and Epidemiology}.
\newblock Springer, Berlin, Heidelberg, 2008.

\bibitem{makris2009}
N.~C. Makris, P.~Ratilal, S.~Jagannathan, Z.~Gong, M.~Andrews, I.~Bertsatos,
  O.~R. God{\o}, R.~W. Nero, and J.~M. Jech.
\newblock Critical population density triggers rapid formation of vast oceanic
  fish shoals.
\newblock {\em Science}, 323(5922):1734--1737, 2009.

\bibitem{Valdinoci2017}
A.~Massaccesi and E.~Valdinoci.
\newblock Is a nonlocal diffusion strategy convenient for biological
  populations in competition?
\newblock {\em J. Math. Biol.}, 74(1-2):113--147, 2017.

\bibitem{matsunaga2017}
Y.~Matsunaga, M.~Noda, H.~Murakawa, K.~Hayashi, A.~Nagasaka, S.~Inoue,
  T.~Miyata, T.~Miura, K.-i. Kubo, and K.~Nakajima.
\newblock Reelin transiently promotes n-cadherin--dependent neuronal adhesion
  during mouse cortical development.
\newblock {\em Proc. Natl. Acad. Sci. USA}, 114(8):2048--2053, 2017.

\bibitem{mayor2016}
R.~Mayor and S.~Etienne-Manneville.
\newblock The front and rear of collective cell migration.
\newblock {\em Nat. Rev. Mol. Cell Biol.}, 17(2):97, 2016.

\bibitem{MFG14}
A.~M. {Middleton}, C.~{Fleck}, and R.~{Grima}.
\newblock {A continuum approximation to an off-lattice individual-cell based
  model of cell migration and adhesion.}
\newblock {\em {J. Theor. Biol.}}, 359:220--232, 2014.

\bibitem{mogilner1999}
A.~Mogilner and L.~Edelstein-Keshet.
\newblock A non-local model for a swarm.
\newblock {\em J. Math. Biol.}, 38(6):534--570, 1999.

\bibitem{MCOe05}
D.~Morale, V.~Capasso, and K.~\"Olschl\"ager.
\newblock An interacting particle system modelling aggregation behavior: from
  individuals to populations.
\newblock {\em J. Math. Biol.}, 50:49--66, 2005.

\bibitem{muntean-toschi}
A.~Muntean and F.~Toschi.
\newblock Collective dynamics from bacteria to crowds: An excursion through
  modeling, analysis and simulation.
\newblock 2014.

\bibitem{15MuTo}
H.~Murakawa and H.~Togashi.
\newblock Continuous models for cell-cell adhesion.
\newblock {\em J. Theor. Biol.}, 374:1--12, 2015.

\bibitem{murray2}
J.~D. Murray.
\newblock {\em Mathematical Biology II: Spatial Models and Biomedical
  Applications}.
\newblock Springer, 3 edition, 2003.

\bibitem{NPT}
G.~Nadin, B.~Perthame, and M.~Tang.
\newblock Can a traveling wave connect two unstable states? {T}he case of the
  nonlocal {F}isher equation.
\newblock {\em C. R. Math. Acad. Sci. Paris}, 349(9-10):553--557, 2011.

\bibitem{okubo2001}
A.~Okubo and S.~A. Levin.
\newblock {\em Diffusion and ecological problems: modern perspectives},
  volume~14.
\newblock Springer Science \& Business Media, 2001.

\bibitem{osborne2017}
J.~M. Osborne, A.~G. Fletcher, J.~M. Pitt-Francis, P.~K. Maini, and D.~J.
  Gavaghan.
\newblock Comparing individual-based approaches to modelling the
  self-organization of multicellular tissues.
\newblock {\em PLoS Comp. Biol.}, 13(2):e1005387, 2017.

\bibitem{othmer-hillen2}
H.~Othmer and T.~Hillen.
\newblock The diffusion limit of transport equations ii: chemotaxis equations.
\newblock {\em SIAM J. Appl. Math.}, 62:1122--1250, 2002.

\bibitem{OuZhang}
C.~Ou and Y.~Zhang.
\newblock Traveling wavefronts of nonlocal reaction-diffusion models for
  adhesion in cell aggregation and cancer invasion.
\newblock {\em Can. Appl. Math. Q.}, 21(1):21--61, 2013.

\bibitem{PAS10}
K.~J. Painter, N.~J. Armstrong, and J.~A. Sherratt.
\newblock The impact of adhesion on cellular invasion processes in cancer and
  development.
\newblock {\em J. Theoret. Biol.}, 264(3):1057--1067, 2010.

\bibitem{painter2015}
K.~J. Painter, J.~M. Bloomfield, J.~A. Sherratt, and A.~Gerisch.
\newblock A nonlocal model for contact attraction and repulsion in
  heterogeneous cell populations.
\newblock {\em Bull. Math. Biol.}, 77(6):1132--1165, 2015.

\bibitem{palsson2000}
E.~Palsson and H.~G. Othmer.
\newblock A model for individual and collective cell movement in dictyostelium
  discoideum.
\newblock {\em Proc. Natl. Acad. Sci. USA}, 97(19):10448--10453, 2000.

\bibitem{Pao92}
C.~V. Pao.
\newblock Blowing-up of solution for a nonlocal reaction-diffusion problem in
  combustion theory.
\newblock {\em J. Math. Anal. Appl.}, 166(2):591--600, 1992.

\bibitem{payne1971}
R.~Payne and D.~Webb.
\newblock Orientation by means of long range acoustic signaling in baleen
  whales.
\newblock {\em Ann. New York Acad. Sci.}, 188:110--141, 1971.

\bibitem{rainey1967}
R.~C. Rainey.
\newblock Radar observations of locust swarms.
\newblock {\em Science}, 157(3784):98--99, 1967.

\bibitem{rejniak2007b}
K.~A. Rejniak and R.~H. Dillon.
\newblock A single cell-based model of the ductal tumour microarchitecture.
\newblock {\em Comp. Math. Meth. Medicine}, 8(1):51--69, 2007.

\bibitem{schaller2005}
G.~Schaller and M.~Meyer-Hermann.
\newblock Multicellular tumor spheroid in an off-lattice voronoi-delaunay cell
  model.
\newblock {\em Phys. Rev. E}, 71(5):051910, 2005.

\bibitem{sekimura1999}
T.~Sekimura, M.~Zhu, J.~Cook, P.~K. Maini, and J.~D. Murray.
\newblock Pattern formation of scale cells in lepidoptera by differential
  origin-dependent cell adhesion.
\newblock {\em Bull. Math. Biol.}, 61(5):807--828, 1999.

\bibitem{SGAP09}
J.~Sherratt, S.~Gourley, N.~Armstrong, and K.~Painter.
\newblock Boundedness of solutions of a non-local reaction-diffusion model for
  adhesion in cell aggregation and cancer invasion.
\newblock {\em Eur. J. Appl. Math.}, 20(1):123--144, 2009.

\bibitem{steinberg1963}
M.~S. Steinberg.
\newblock Reconstruction of tissues by dissociated cells.
\newblock {\em Science}, 141(3579):401--408, 1963.

\bibitem{steinberg2007}
M.~S. Steinberg.
\newblock Differential adhesion in morphogenesis: a modern view.
\newblock {\em Curr. Opin. Gen. \& Devel.}, 17(4):281--286, 2007.

\bibitem{sumpter2010}
D.~J.~T. Sumpter.
\newblock {\em Collective animal behavior}.
\newblock Princeton University Press, 2010.

\bibitem{taylor2017}
H.~B. Taylor, A.~Khuong, Z.~Wu, Q.~Xu, R.~Morley, L.~Gregory, A.~Poliakov,
  W.~R. Taylor, and D.~G. Wilkinson.
\newblock Cell segregation and border sharpening by eph
  receptor--ephrin-mediated heterotypic repulsion.
\newblock {\em J. Roy. Soc. Interface}, 14(132):20170338, 2017.

\bibitem{topaz2006}
C.~M. Topaz, A.~L. Bertozzi, and M.~A. Lewis.
\newblock A nonlocal continuum model for biological aggregation.
\newblock {\em Bull. Math. Biol.}, 68(7):1601, 2006.

\bibitem{townes1955}
P.~L. Townes and J.~Holtfreter.
\newblock Directed movements and selective adhesion of embryonic amphibian
  cells.
\newblock {\em J. Exper. Zool.}, 128(1):53--120, 1955.

\bibitem{trush2019}
O.~Trush, C.~Liu, X.~Han, Y.~Nakai, R.~Takayama, H.~Murakawa, J.~A. Carrillo,
  H.~Takechi, S.~Hakeda-Suzuki, T.~Suzuki, et~al.
\newblock N-cadherin orchestrates self-organization of neurons within a
  columnar unit in the drosophila medulla.
\newblock {\em J. Neurosci.}, 39(30):5861--5880, 2019.

\bibitem{turing1952}
A.~M. Turing.
\newblock The chemical basis of morphogenesis.
\newblock {\em Philos. Trans. Roy. Soc. London Ser. B}, 237:37--72, 1952.

\bibitem{vol2}
V.~Volpert.
\newblock {\em Elliptic partial differential equations. {V}ol. 2}, volume 104
  of {\em Monographs in Mathematics}.
\newblock Birkh\"{a}user/Springer Basel AG, Basel, 2014.
\newblock Reaction-diffusion equations.

\bibitem{wang-collective}
X.~Wang, A.~Enomoto, N.~Asai, T.~Kato, and M.~Takahashi.
\newblock Collective invasion of cancer: Perspectives from pathology and
  development.
\newblock {\em Pathol. Int.}, 66(4):183--192, 2016.

\bibitem{Xiang}
T.~Xiang.
\newblock A study on the positive nonconstant steady states of nonlocal
  chemotaxis systems.
\newblock {\em Discrete Contin. Dyn. Syst. Ser. B}, 18(9):2457--2485, 2013.

\bibitem{ZhigunKS}
A.~Zhigun.
\newblock Generalised global supersolutions with mass control for systems with
  taxis.
\newblock {\em SIAM J. Math. Anal.}, 51(3):2425--2443, 2019.

\end{thebibliography}

\end{document}